\documentclass[11pt,reqno]{amsart}
\usepackage{indentfirst, amssymb,amsmath,amsthm, mathrsfs,cite}
\newtheorem*{theoA}{Theorem A}
\newtheorem*{theoB}{Theorem B}
\newtheorem*{theoC}{Theorem C}
\newtheorem*{theoD}{Theorem D}

\newtheorem{theo}{Theorem}[section]
\newtheorem{lem}{Lemma}[section]
\newtheorem{cor}{Corollary}[section]
\newtheorem{remark}{Remark}[section]
\newtheorem{ex}{Example}[section]

\newtheorem*{propA}{Proposition A}
\newtheorem*{propB}{Proposition B}
\newtheorem*{propC}{Proposition C}

\newtheorem{open problem}{Open problem}[section]
\newcommand{\pa}{\partial}
\newcommand{\ol}{\overline}
\newcommand{\be}{\begin{equation}}
\newcommand{\ee}{\end{equation}}
\newcommand{\bs}{\begin{small}}
\newcommand{\es}{\end{small}}
\newcommand{\beas}{\begin{eqnarray*}}
\newcommand{\eeas}{\end{eqnarray*}}
\newcommand{\bea}{\begin{eqnarray}}
\newcommand{\eea}{\end{eqnarray}}
\renewcommand{\epsilon}{\varepsilon}
\numberwithin{equation}{section}

\begin{document}
\title[Solutions for certain Fermat-type PDDE\lowercase{s}]{Solutions for certain Fermat-type PDDE\lowercase{s} concerning an open problem of Xu and Wang}
\author[H. Y. XU, R. Mandal and R. Biswas]{Hong Yan Xu, Rajib Mandal and Raju Biswas}
\date{}
\address{School of Arts and Sciences, Suqian University, Suqian, Jiangsu 223800, China.}
\email{xhyhhh@126.com}
\address{Department of Mathematics, Raiganj University, Raiganj, West Bengal-733134, India.}
\email{rajibmathresearch@gmail.com}
\address{Department of Mathematics, Raiganj University, Raiganj, West Bengal-733134, India.}
\email{rajubiswasjanu02@gmail.com}
\maketitle
\let\thefootnote\relax
\footnotetext{Type set by \AmS -\LaTeX}
\begin{abstract} 
The objective of this study is to ascertain the existence and forms of the finite order meromorphic and entire functions of several complex variables satisfying some certain Fermat-type partial differential-difference equations by considering the more general forms of the PDDEs in an open problem on $\mathbb{C}^2$ due to Xu and Wang (Notes on the existence of entire solutions for several partial differential-difference equations, Bull. Iran. Math. Soc., 47, 1477-1489 (2020)). We provide examples to illustrate the results.
\end{abstract}
\section{Introduction, Definitions and Results}
By a meromorphic function $f$ on $\mathbb{C}^n$ ($n\in\mathbb{N}$), we mean that $f$ can be written as a quotient of two holomorphic functions without common zero sets in 
$\mathbb{C}^n$.  Notationally, we write $f:=\frac{g}{h}$, where $g$ and $h$ are relatively prime holomorphic functions on $\mathbb{C}^n$ such that $h \not\equiv 0$ and 
$f^{-1}(\infty)\not=\mathbb{C}^n$. In particular, the entire function of several complex variables are holomorphic throughout $\mathbb{C}^n$. 

Let $z=(z_1,z_2,\ldots,z_n)\in\mathbb{C}^n$, $a\in\mathbb{C}\cup\{\infty\}$, $k\in\mathbb{N}$ and $r>0$. We consider some notations from \cite{14,S1,Y1}. Let $\ol B_n(r):=\{z\in\mathbb{C}^n: |z|\leq r\}$, where $|z|^2:=\sum_{j=1}^n|z_j|^2$. 
The exterior derivative splits $d:=\pa+\ol{\pa}$ and twists to $d^c:=\frac{i}{4\pi}(\ol\pa-\pa)$. The standard Kaehler metric on $\mathbb{C}^n$ is 
given by $v_n(z):=dd^c|z|^2$. Define $\omega_n(z):=dd^c\log |z|^2\geq 0$ and $\sigma_n(z):=d^c \log |z|^2\wedge \omega_n^{n-1}(z)$ on $\mathbb{C}^n\setminus\{0\}$. Thus $\sigma_n(z)$ defines a positive measure on $\pa B_n:=\{z\in\mathbb{C}^n: |z|= r\}$ with total measure $1$.
The zero-multiplicity of a holomorphic function $h$ at a point 
$z\in \mathbb{C}^n$ is defined to be the order of vanishing of $h$ at $z$ and denoted by $\mathcal{D}_h^0(z)$. A divisor of $f$ on $\mathbb{C}^n$ is an integer valued function which is locally the difference between the 
zero-multiplicity functions of $g$ and $h$ and it is denoted by $\mathcal{D}_f:=\mathcal{D}_g^0-\mathcal{D}_h^0$ (see, P. 381, \cite{51}). Let $a\in\mathbb{C}\cup\{\infty\}$ be such that $f^{-1}(a)\not=\mathbb{C}^n$. Then the $a$-divisor $\nu_f^a$ of $f$ 
is the divisor associated with the holomorphic functions $g-ah$ and $h$ (see, P. 346, \cite{14} and P. 12, \cite{12}). 
Ye \cite{Y1} has defined the counting function and the valence function with respect to $a$ respectively as follows:
\beas n(r,a,f):=r^{2-2n}\int_{S(r)} \nu_f^a v_n^{n-1}\;\text{and}\;N(r,a,f):=\int_0^r\frac{ n(r,a,f)}{t}dt. \eeas
We write 
\beas N(r,a,f)=\left\{\begin{array}{ll}
&N\left(r,\frac{1}{f-a}\right),\;\text{when}\;a\not=\infty\\[2mm]
&N(r,f),\;\text{when}\;a=\infty.\end{array}\right.\eeas 
The proximity function \cite{14,Y1} of $f$ is defined as follows :
\beas\begin{array}{ll}
m(r,f):=\int_{\pa B_n(r)} \log^+|f(z)|\sigma_n(z),\;\text{when}\;a=\infty\\
m\left(r, \frac{1}{f-a}\right):=\int_{\pa B_n(r)} \log^+\frac{1}{|f(z)-a|}\sigma_n(z),\text{when}\;a\not=\infty.\end{array}\eeas

The Nevanlinna characteristic function is defined by $T(r,f)=N(r,f)+m(r,f)$, which is increasing for $r$. The order of a meromorphic function $f$ is denoted by $\rho(f)$ and is defined by 
\beas \rho(f)=\varlimsup\limits_{r\to \infty}\frac{\log^+ T(r,f)}{\log r},\;\;\text{where}\;\log^{+}x=\max\{\log x,0\}.\eeas
The exceptional sets are throughout needed in the Nevanlinna theory. Typically, it means considering the linear measure $m(E):=\int_E dt$ and the logarithmic measure 
$l(E):=\int_{E\cap [1, \infty)} dt/t$ for a set $E\subset [0,\infty)$. 
Recall that a meromorphic function $\alpha$ is said to be a small function of $f$, if $T(r,\alpha)=S(r,f)$, where $S(r,f)$ is any quantity that satisfies 
$S(r,f)=o(T(r,f))$ as $r\rightarrow \infty$, possibly outside of a set of $r$ of finite linear measure. 
For further details, we refer to \cite{51,12i,12,15,19,21,S1,Y1} and the references therein.
Given a meromorphic function $f(z)$ on $\mathbb{C}^n$, $f(z+c)$ is called a shift of $f$ and $\Delta (f)=f(z+c)-f(z)$ is called a difference operator of $f$, where $c(\not=0)\in\mathbb{C}^n$.\\[2mm]
\indent An equation is called a partial differential equation (in brief, PDE) if the equation contains partial derivatives of $f$ whereas if the equation also contains shifts or differences of $f$, then the equation is called a partial differential-difference equation (in brief, PDDE).
We now consider the Fermat-type equation
\bea\label{eq1.1} f^n(z)+g^n(z)=1,\;\text{where}\;n\in\mathbb{N}.\eea 
We summarize the classical results for solutions of the equation (\ref{eq1.1}) on $\mathbb{C}$ in the following:
\begin{propA}
(i)\cite{2201,2001} The equation (\ref{eq1.1}) with $n=2$ has the non-constant entire solutions $f(z)=\cos(\eta(z))$ and $g(z)=\sin(\eta(z))$, where $\eta(z)$ is any entire function. No other solutions exist.\\
(ii)\cite{2201,7,17} For $n\geq 3$, there are no non-constant entire solutions of (\ref{eq1.1}) on $\mathbb{C}$.
\end{propA}
\begin{propB}
(i)\cite{2201} The equation (\ref{eq1.1}) with $n=2$ has 
the non-constant meromorphic solutions 
$f=\frac{2\omega}{1+\omega^2}$ and $g=\frac{1-\omega^2}{1+\omega^2}$, where $\omega$ is an arbitrary meromorphic function on $\mathbb{C}$.\\ 
(ii)\cite{3601,7} The equation (\ref{eq1.1}) with $n=3$ has
the non-constant meromorphic solutions 
$f=\frac{1}{2\wp(h)}\left(1+\frac{\wp'(h)}{\sqrt{3}}\right)$, $g=\frac{1}{2\wp(h)}\left(1-\frac{\wp'(h)}{\sqrt{3}}\right)$, 
where $\wp(z)$ denotes the Weierstrass elliptic $\wp$-function with periods $\omega_1$ and $\omega_2$ is defined as
\beas \wp \left(z;\omega_1,\omega_2\right)=\frac{1}{z^2}+\sum\limits_{{\mu,\nu;\mu ^2+\nu ^2\not=0}} \left\{ \frac{1}{\left(z+\mu\omega_1+\nu\omega_2\right)^2}-\frac{1}{\left(\mu\omega_1+\nu\omega_2\right)^2}\right\},\eeas
which is even and satisfying, after appropriately choosing $\omega_1$ and $\omega_2$, $(\wp')^2=4\wp^3-1$.\\
(iii)\cite{2201,7} For $n\geq 4$, there are no non-constant meromorphic solutions of (\ref{eq1.1}) on $\mathbb{C}$.
\end{propB}
Numerous researchers have shown their interest to investigate on the Fermat-type equations for entire and meromorphic solutions from last two decades by taking some variation of (\ref{eq1.1}). Yang and Li \cite{170} was the pioneer for introducing the study on transcendental meromorphic solutions of Fermat-type differential equation on $\mathbb{C}$. Liu \cite{Liu1} was the first who investigated on meromorphic solutions of Fermat-type difference equation as well as differential-difference equations on $\mathbb{C}$. For the leading and recent developments in these directions, we refer to the reader to \cite{Raj1,Raj2,601,2201,7,28,16,Liu2,602,603,604,600,100} and the references therein.

The basic conclusions \textrm{Propositions A} and \textrm{B} of the Fermat-type equation (\ref{eq1.1}) on $\mathbb{C}$ were also extended to the case of several complex variables and the following is the summarization.
\begin{propC} \cite[\textrm{Theorem 2.3}]{700} \cite[\textrm{Theorem 1.3}]{20} Let $h:\mathbb{C}^n\to\mathbb{C}$ be a non-constant entire function and $\mathbb{P}^1=\mathbb{C}\cup\{\infty\}$. Then the non-constant entire and meromorphic solutions of the equation (\ref{eq1.1}) on $\mathbb{C}^n$ are characterized as follows:\\
(i) when $m = 2$, the entire solutions are $f = cos(h)$ and $g = sin(h)$;\\
(ii) when $m > 2$, there are no non-constant entire solutions;\\
(iii) when $m=2$, the meromorphic solutions are $f=\frac{2\omega}{1+\omega^2}$ and $g=\frac{1-\omega^2}{1+\omega^2}$, where $\omega:\mathbb{C}^n\to \mathbb{P}^1$ is a non-constant meromorphic function;\\ 
(iv) when $m=3$, the meromorphic solutions are
$f=\frac{1}{2\wp(h)}\left(1+\frac{\wp'(h)}{\sqrt{3}}\right)$ and $g=\frac{1}{2\wp(h)}\left(1-\frac{\wp'(h)}{\sqrt{3}}\right)$, where $\wp(z)$ denotes the Weierstrass elliptic $\wp$-function satisfying the relation $(\wp')^2=4\wp^3-1$. Note that $\wp:\mathbb{C}\to \mathbb{P}^1$ so that $\wp\circ h:\mathbb{C}^n\to \mathbb{P}^1$, i.e., $f:\mathbb{C}^n\to \mathbb{P}^1$;\\
(v) when $m> 3$, there are no non-constant meromorphic solutions.
\end{propC}
Now researchers have been focusing their attention to investigate on the Fermat-type PDDEs for entire and meromorphic solutions. Let \bea\label{esl}\sum_{i=1}^n\left(\frac{\pa u}{\pa z_i}\right)^m=1\eea
be the certain non-linear first order PDE introducing from the analogy with the Fermat-type equation $\sum_{i=1}^n\left(f_i\right)^m=1$, where $u:\mathbb{C}^n\to \mathbb{C}$, $z_i\in\mathbb{C}$, $f_i:\mathbb{C}\to \mathbb{C}$, and $m,n\geq 2$.
In 1999, Saleeby \cite{27} first started to study about the solutions of the Fermat-type PDEs and obtained the results for entire solutions of (\ref{esl}) on $\mathbb{C}^2$. Afterwards, in 2004, Li \cite{28} extended these results to $\mathbb{C}^n$.\\
\indent In 2008, Li \cite{16} considered the equation (\ref{eq1.1}) with $n=2$ and showed that meromorphic solutions $f$ and $g$ of that equation on 
$\mathbb{C}^2$ must be constant if and only if $\pa f/\pa z_2$ and $\pa g/\pa z_1$ have the same zeros (counting multiplicities). 
If $f=\pa u/\pa z_1$ and $g=\pa u/\pa z_2$, then any entire solutions of the partial differential equations $\left(\pa u/\pa z_1\right)^2+\left(\pa u/\pa z_2\right)^2=1$ on $\mathbb{C}^2$ are necessarily linear \cite{13}.\par
In 2018, Xu and Cao \cite{23,26} was the first who considered both difference operators and differential operators in Fermat-type equations of two complex variables and obtained the following results.
\begin{theoA}\cite{23,26} The PDDE
\bea\label{xc} \left(\frac{\pa f(z)}{\pa z_1}\right)^n+f^m(z+c)=1,\;\text{where}\; c=(c_1,c_2)\in\mathbb{C}^2,\eea
doesn't have any finite order transcendental entire solution of two complex variables $z_1$ and $z_2$, where $m,n\in\mathbb{N}$ are distinct.\end{theoA}
\begin{theoB}\cite{23,26} Let $m=n=2$. Then any transcendental entire solution with finite order of (\ref{xc}) must have the form $f\left(z_1,z_2\right)=\sin (Az_1+Bz_2+H(z_2))$, where $A,B\in\mathbb{C}$ satisfying $A^2=1$ and $Ae^{i(Ac_1+Bc_2)}=1$, and $H\left(z_2\right)$ is a polynomial in one variable $z_2$ such that $H\left(z_2\right)\equiv H\left(z_2+c_2\right)$. In a 
special case, if $c_2\not=0$, then $f\left(z_1,z_2\right)=\sin (Az_1+Bz_2+ C)$, where $C\in\mathbb{C}$.\end{theoB}
The authors \cite{23,26} also proved that, if $c_1=c_2=0$ in PDDE (\ref{xc}) with $m=n=2$, then any finite order transcendental entire solution of (\ref{xc}) is of the form $f\left(z_1,z_2\right)=\sin \left(z_1+g\left(z_2\right)\right)$, where $g\left(z_2\right)$ is a polynomial in one variable $z_2$.

The authors \cite{23,26} also obtained the first result on the meromorphic solutions of (\ref{xc}) and it described as follows.
\begin{theoC}\cite{23,26} Let $m=n=2$ and $c=(c_1,c_2)\in\mathbb{C}^2$. Then any non-constant meromorphic solution of (\ref{xc}) must have the form $f(z)=\frac{h(z-c)-\frac{1}{h(z-c)}}{2i}$, where $h$ is a non-zero meromorphic function on $\mathbb{C}^2$ satisfying $i\left(h(z+c)+\frac{1}{h(z+c)}\right)=\frac{\pa h(z)}{\pa z_1}\left(1+\frac{1}{h(z)}\right)$.
In a special case, where $c_1=c_2=0$, we have $f(z)=\sin \left(z_1-ia(z_2)\right)$, where $a(z_2)$ is a meromorphic function in one complex variable $z_2$.\end{theoC}
In 2020, Xu and Wang \cite{24} took some variations of the equation (\ref{xc}), replacing $\frac{\pa f(z)}{\pa z_1}$ by the term $\frac{\pa f(z)}{\pa z_1}+\frac{\pa f(z)}{\pa z_2}$. Actually they considered the following PDDEs on $\mathbb{C}^2$:    
\bea\label{xw} \left(\frac{\pa f(z)}{\pa z_1}+\frac{\pa f(z)}{\pa z_2}\right)^n+f^m(z+c)=1,\;\text{where}\;c=(c_1,c_2)\in\mathbb{C}^2, \eea
and proved the following results.
\begin{theoD}\cite{24} Let $m,n\in\mathbb{N}$ be distinct. 
Then (\ref{xw}) does not have any finite order transcendental entire solution of two complex variables $z_1$ and $z_2$, whenever $m>n$ or $n>m\geq 2$.\end{theoD}
Moreover, the authors \cite{24} obtained the first result on finite order entire solutions of two complex variables, where the combination were $n=2$, $m=1$ in (\ref{xc}) and (\ref{xw}), i.e.,
\begin{small}\bea\label{equ1}&&\left(\frac{\pa f(z)}{\pa z_1}\right)^2+f(z+c)=1\\  
\text{and} &&\label{equ2}\left(\frac{\pa f(z)}{\pa z_1}+\frac{\pa f(z)}{\pa z_2}\right)^2+f(z+c)=1\eea \end{small}
have the finite order transcendental entire solutions respectively
$f(z_1,z_2)=1-\frac{1}{4}c_1^2-\frac{1}{4}z_1^2+\frac{c_1}{2c_2}z_1z_2-\frac{c_1^2}{2c_2}(z_2-c_2)+(z_1-c_1)G_1(z_2)-\left[\frac{c_1}{2c_2}\left(z_2-c_2\right)+G_1(z_2)\right]^2$ and 
$f(z_1,z_2)=1-\frac{1}{4}c_1^2-\frac{1}{4}z_1^2+z_1\left[G_2\left(z_2-z_1\right)+a_3\left(z_2-z_1\right)\right]-c_1G_2\left(z_2-z_1\right)-a_3c_1\left[z_2-z_1-\left(c_2-c_1\right)\right]-\left[G_2\left(z_2-z_1\right)+a_3\left(z_2-z_1-\left(c_2-c_1\right)\right)\right]^2$, where $G_1(z_2)$, $G_2(z_2-z_1)$ are finite order transcendental entire period functions with period $c_2$, $c_2-c_1$ respectively and $a_3=\frac{c_1}{2\left(c_2-c_1\right)}$.\\
Lastly, Xu and Wang \cite{24} posed the following open problem in their paper.
\begin{open problem}\label{op1} Whether there exists the finite order transcendental entire solutions of two complex variables $z_1$ and $z_2$ for the equations (\ref{xc}) and (\ref{xw}) in the case $n>2$ and $m=1$ or not ?\end{open problem}
As far as we know, this open problem is not solved till now. Our main aim of this paper is to solve the  \textrm{Open problem \ref{op1}}. In this paper, we consider the more compact forms of (\ref{xc}) and (\ref{xw}), and then solve these equations for the finite order transcendental entire functions of several complex variables. Thus the \textrm{Open problem \ref{op1}} has been solved in this paper.
\section{The main results}
Let $I=(i_1,i_2,\ldots,i_k)\in\mathbb{Z}^k_+$ be a multi-index with length $\Vert I\Vert=\sum_{j=1}^k i_j$, $\mathbb{Z}_+=\mathbb{N}\cup\{0\}$ and $\pa^I f=\frac{\pa^{\Vert I\Vert} f}{\pa z_1^{i_1}\cdots\pa z_k^{i_k}}$.
Now, any polynomial $\mathcal{Q}(z)$ of several complex variables of degree $d$ can 
be expressed as $\mathcal{Q}(z)=\sum_{\Vert I\Vert=0}^d a_I z_1^{i_1}\cdots z_n^{i_n}$, where $a_{I}\in\mathbb{C}$ such that $a_{I}$ are not all zero at a time for $\Vert I\Vert=d$.  
Let $G(z)$ denotes the partial differential function of finite order transcendental meromorphic $f$ function on $\mathbb{C}^n$ with $N(r,f)=S(r,f)$ involving $n$($\in\mathbb{N}$) different homogeneous terms on $\mathbb{C}^n$ such that
\bs\bea\label{re1}G(z)&=&\sum\limits_{m=1}^{n}\sum\limits_{\Vert I\Vert=m}a_{I}(z)\pa^I f(z)\\
&=&\left(a_{(n,0\ldots,0)}(z)\frac{\pa^nf(z)}{\pa z_1^n}+\cdots+a_{(1,1,1,\ldots,1)}(z)\frac{\pa^nf(z)}{\pa z_1\pa z_2\cdots\pa z_n}+\cdots+a_{(0,\ldots,0,n)}(z)\frac{\pa^nf(z)}{\pa z_n^n}\right)\nonumber\\
&&+\left(a_{(n-1,0,\ldots,0)}(z)\frac{\pa^{n-1}f(z)}{\pa z_1^{n-1}}
+\cdots+a_{(0,\ldots,0,n-1)}(z)\frac{\pa^{n-1}f(z)}{\pa z_n^{n-1}}\right)\nonumber\\
&&+\cdots+\left(a_{(1,0,\ldots,0)}(z)\frac{\pa f(z)}{\pa z_1}+a_{(0,1,\ldots,0)}(z)\frac{\pa f(z)}{\pa z_2}+\cdots+a_{(0,\ldots,0,1)}(z)\frac{\pa f(z)}{\pa z_n}\right),\nonumber\eea\es
where $z=\left(z_1,z_2,\ldots,z_n\right)$ and $a_{I}(z)$ are small functions of $f(z)$ of several complex variables such that $a_{I}(z)$ are not all identically zero at a time.
We now investigate about the existence of solutions of the following Fermat-type PDDE on $\mathbb{C}^n$:
\bea\label{fg} G^{m_1}(z)+\alpha(z)(\Delta (f))^{m_2}=\beta(z),\eea
where $m_1,m_2\in\mathbb{N}$, $c=\left(c_1,c_2,\ldots,c_n\right)\in\mathbb{C}^n$ with $c\not=0$ and $\alpha(z)(\not\equiv 0)$, $\beta(z)(\not\equiv 0)$ are small functions of $f$ of several complex variables. 
For existence of solution of (\ref{fg}), we obtain the following result.
\begin{theo}\label{th} Let $m_1,m_2\in\mathbb{N}$ be distinct. Then (\ref{fg}) does not have any finite order transcendental meromorphic solution $f$ of several complex 
variables, where $N(r,f)=S(r,f)$ and $m_2>m_1$ or $m_1>m_2\geq 2$.
\end{theo}
Clearly \textrm{Theorem \ref{th}} improves as well as generalizes significantly both \textrm{Theorems A} and \textrm{D}.\par
Let $m_1\in\mathbb{N}\setminus\{1\}$. Corresponding to the \textrm{Open problem \ref{op1}}, we now consider the entire functions of several complex variables with finte order satisfying the following Fermat-type PDDEs
\bea\label{fte}&&\left(\frac{\pa f(z)}{\pa z_1}\right)^{m_1}+\Delta (f)=\varphi(z_2,z_3,\ldots,z_n)\\\text{and}
\label{ftee}&&\left(\frac{\pa f(z)}{\pa z_1}+\frac{\pa f(z)}{\pa z_2}\right)^{m_1}+\Delta (f)=\varphi(z_3,z_4,\ldots,z_n),\eea
where $c(\not=0)\in\mathbb{C}^n$, $\varphi(z_2,z_3,\ldots,z_n)(\not\equiv 0)$ and $\varphi(z_3,z_4,\ldots,z_n)(\not\equiv 0)$ are finite order entire functions. 

For the finite order transcendental entire functions of several complex variables satisfying (\ref{fte}) and (\ref{ftee}), we obtain the following results respectively.
\begin{theo}\label{T1} Let $f$ be a finite order transcendental entire function on $\mathbb{C}^n$ that satisfies (\ref{fte}). If $m_1=2$, then the entire solution of (\ref{fte}) has one of the following form:
\item[(I)]
\bs\beas f(z)=\varphi(z_2-c_2,z_3-c_3,\ldots,z_n-c_n)-\left(-\frac{1}{2}(z_1-c_1)+g_1(z_2-c_2,z_3-c_3,\ldots,z_n-c_3)\right)^2,\eeas\es
where $g_1(z_2,z_3,\ldots,z_n)$ is a polynomial such that $g_1(z_2+c_2,z_3+c_3,\ldots,z_n+c_n)\equiv g_1(z_2,z_3,\ldots,z_n)+\frac{c_1}{2}$ holds and $\varphi(z_2,z_3,\ldots,z_n)$ is a finite order transcendental entire function in $z_2,z_3,\ldots,z_n$;
\item[(II)]  
\beas&& f(z)=(z_1-c_1)\left(g_2(z_2,z_3,\ldots,z_n)+\frac{c_1}{2\tau}\omega\right)-\left(g_2(z_2,z_3,\ldots,z_n)+\frac{c_1}{2\tau}(\omega-\tau)\right)^2\\&&+\frac{1}{4}\left(c_1^2-z_1^2\right)+\varphi(z_2-c_2,z_3-c_3,\ldots,z_n-c_n),\eeas
where $g_2(z_2,z_3,\ldots,z_n)$ is a finite order transcendental entire periodic function in $z_2,z_3,\\\ldots,z_n$ with period $(c_2,c_3,\ldots,c_n)\in\mathbb{C}^{n-1}$, $\varphi(z_2,z_3,\ldots,z_n)$ is a finite order entire function, $\omega=\sum_{j=2}^nz_j$ and $\tau=\sum_{j=2}^n c_j\not=0$.\\
If $m_1\geq 3$, then the equation (\ref{fte}) does not have any finite order transcendental entire solution.
\end{theo}
In particular, if $\varphi(z_2,z_3,\ldots,z_n)\equiv 1$, then we obtain the following corollary.
\begin{cor}\label{cor1} Let $f$ be a finite order transcendental entire function on $\mathbb{C}^n$ that satisfies (\ref{fte}) with $\varphi(z_2,z_3,\ldots,z_n)\equiv 1$. If $m_1=2$, then the equation (\ref{fte}) has the entire solution of the form
\begin{small}\beas&& f(z)=1+\frac{1}{4}\left(c_1^2-z_1^2\right)+\frac{c_1}{2\tau}z_1\omega+z_1g_2(z_2,z_3,\ldots,z_n)-\left(g_2(z_2,z_3,\ldots,z_n)+\frac{c_1}{2\tau}(\omega-\tau)\right)^2\\&&-c_1\left(g_2(z_2,z_3,\ldots,z_n)+\frac{c_1}{2\tau}\omega\right),\eeas \end{small}
where $g_2(z_2,z_3,\ldots,z_n)$ is a finite order transcendental entire periodic function in $z_2,z_3,\ldots,z_n$ with period $(c_2,c_3,\ldots,c_n)\in\mathbb{C}^{n-1}$, $\omega=\sum_{j=2}^nz_j$ and $\tau=\sum_{j=2}^nc_j\not=0$.\\
If $m_1\geq 3$, then the equation (\ref{fte}) does not have any finite order transcendental entire solution.
\end{cor}
The following examples related to \textrm{Theorem \ref{T1}} are reasonable.
\begin{ex} Let 
\bs\beas &&f(z_1,z_2,\ldots,z_5)=\pi i+z_3-z_4+z_5+e^{z_2+z_3-2z_4}-\frac{1}{4}(\pi^2+z_1^2)+(z_1-\pi i)e^{5z_2z_3-2z_2z_4+z_5+9}\\
&&+\frac{1}{18}(z_1-\pi i)(z_2+z_3+z_4+z_5)-\left(e^{5z_2z_3-2z_2z_4+z_5+9}+\frac{1}{18}(z_2+z_3+z_4+z_5-9\pi i)\right)^2.\eeas\es 
It is easy to see that $f$ is a transcendental entire function on $\mathbb{C}^5$ with $\rho(f)=2$ and satisfying the equation \beas \left(\frac{\partial f(z_1,z_2,\dots,z_5)}{\partial z_1}\right)^2+f(z_1+c_1,z_2+c_2,\ldots,z_5+c_5)=e^{z_2+z_3-2z_4}+z_3-z_4+z_5,\eeas
where $c=(\pi i,0,2\pi i,5\pi i,2\pi i)$.\end{ex}
\begin{ex} Let $f(z_1,z_2,\ldots,z_5)=\frac{4+\pi^2}{4}-\frac{1}{4}z_1^2+(z_1-\pi)e^{7z_2-2z_3+5z_4-3z_5+1}+\frac{1}{3i}(z_1-\pi)(z_2+z_3+z_4+z_5)-\left(e^{7z_2-2z_3+5z_4-3z_5+1}+\frac{1}{3i}(z_2+z_3+z_4+z_5-3\pi i/2)\right)^2$ be a transcendental entire function on $\mathbb{C}^5$. Then $\rho(f)=1$ and clearly $f$ satisfies the equation 
\beas \left(\frac{\partial f(z_1,z_2,\ldots,z_5)}{\partial z_1}\right)^2+f(z_1+c_1,z_2+c_2,\ldots,z_5+c_5)=1,\eeas 
where $c=(\pi,\pi,\pi i/2,-\pi i, \pi i)$.\end{ex}
\begin{ex} Consider $f(z_1,z_2,z_3,z_4)=e^{z_2+2z_3-z_4-2}-\left(-z_1/2-z_2+z_3+z_4\right)^2$. Then $f$ is of order $1$ and satisfies the equation 
\beas \left(\frac{\partial f(z_1,z_2,z_3,z_4)}{\partial z_1}\right)^2+f(z_1+c_1,z_2+c_2,z_3+c_3,z_4+c_4)=e^{z_2+2z_3-z_4},\eeas
where $c=(14 ,1,3,5)$.\end{ex}
\begin{ex} Consider $f(z_1,z_2,z_3)=1-\frac{1}{4}z_1^2+z_1e^{z_2+z_3}-e^{2z_2+2z_3}$. Then $f$ is of order $1$ and satisfies the equation 
\beas \left(\frac{\partial f(z_1,z_2,z_3)}{\partial z_1}\right)^2+f(z_1+c_1,z_2+c_2,z_3+c_3)=1,\;\text{where}\;c=(0,\pi i,\pi i).\eeas\end{ex}
\begin{theo}\label{T2}Let $f$ be a finite order transcendental entire function on $\mathbb{C}^n$ that satisfies satisfies (\ref{ftee}). If $m_1=2$, then the entire solution of (\ref{ftee}) has one of the following form:
\item[(I)] \beas &&f(z)=-\left(-\frac{1}{2}(z_1-c_1)+g_2(z_2-z_1-c_2+c_1,z_3-c_3,\cdots,z_n-c_n)\right)^2\\&&+\varphi(z_3-c_3,z_4-c_4,\ldots,z_n-c_n),\eeas
where $g_2$ is a polynomial satisfies $g_2(z_2-z_1+c_2-c_1,z_3+c_3,\cdots,z_n+c_n)\equiv g_2(z_2-z_1,z_3,\cdots,z_n)+\frac{c_1}{2}$ with $\frac{\pa g_2}{\pa z_1}+\frac{\pa g_2}{\pa z_2}\equiv 0$ and $\varphi(z_3,z_4,\ldots,z_n)$ is a finite order transcendental entire function;
\item[(II)] \beas &&f(z)=\frac{1}{4}\left(c_1^2-z_1^2\right)+(z_1-c_1)\left(g_4(z_2-z_1,z_3,\cdots,z_n)+\frac{c_1}{2\tau}(z_2-z_1+z_3+\cdots+z_n)\right)\\&&-\left(g_4(z_2-z_1,z_3,\cdots,z_n)+\frac{c_1}{2\tau}(z_2-z_1+z_3+\cdots+z_n-\tau)\right)^2\\&&+\varphi(z_3-c_3,z_4-c_4,\ldots,z_n-c_n),\eeas
where $g_4(z_2-z_1,z_3,\cdots,z_n)$ is a finite order transcendental entire function with period \\$(c_2-c_1,c_3,\ldots,c_n)$ with $\frac{\pa g_4}{\pa z_1}+\frac{\pa g_4}{\pa z_2}\equiv0$, $\varphi(z_3,z_4,\ldots,z_n)$ is a finite order entire function, $\omega=z_2-z_1+z_3+\cdots+z_n$, $\tau=c_2-c_1+c_3+c_4+\cdots+c_n\not=0$.\\
If $m_1\geq 3$, then the equation (\ref{ftee}) does not have any finite order transcendental entire solution.
\end{theo}
In particular, if $\varphi(z_3,z_4,\ldots,z_n)\equiv 1$, then we obtain the following corollary.
\begin{cor}\label{cor2} Let $f$ be a finite order transcendental entire function on $\mathbb{C}^n$ that satisfies (\ref{ftee}) with $\varphi(z_3,z_4,\ldots,z_n)\equiv 1$. If $m_1=2$, then the equation (\ref{ftee}) has the entire solution of the form 
\begin{small}\beas &&f(z)=1+z_1\left(g_4(z_2-z_1,z_3,\cdots,z_n)+\frac{c_1}{2\tau}\omega\right)-c_1\left(g_4(z_2-z_1,z_3,\cdots,z_n)+\frac{c_1}{2\tau}(\omega-\tau)\right)\\&&-\left(g_4(z_2-z_1,z_3,\cdots,z_n)+\frac{c_1}{2\tau}(\omega-\tau)\right)^2-\frac{1}{4}\left(c_1^2+z_1^2\right),\eeas\end{small}
where $\omega=z_2-z_1+z_3+\cdots+z_n$, $\tau=c_2-c_1+c_3+c_4+\cdots+c_n\not=0$ and $g_4(z_2-z_1,z_3,\cdots,z_n)$ is a finite order transcendental entire function with period $(c_2-c_1,c_3,\ldots,c_n)\in\mathbb{C}^{n-1}$ satisfying $\frac{\pa g_4}{\pa z_1}+\frac{\pa g_4}{\pa z_2}\equiv 0$.\\
If $m_1\geq 3$, then the equation (\ref{ftee}) does not have any finite order transcendental entire solution.
\end{cor}
The following examples related to \textrm{Theorem \ref{T2}} are reasonable.
\begin{ex} Consider 
$f(z_1,z_2,z_3,z_4)=(z_3-2)e^{2z_3+z_4-8}-\left(9z_1/2-5z_2+7z_3-2z_4\right)^2$. Then $f$ is of order $1$ on $\mathbb{C}^4$ and satisfies the equation 
\beas \left(\frac{\partial f(z_1,z_2,z_3,z_4)}{\partial z_1}\right)^2+f(z_1+c_1,z_2+c_2,z_3+c_3,z_4+c_4)=z_3e^{2z_3+z_4},\eeas
where $c=(2,3,2,4)$.\end{ex}
\begin{ex} Clearly $f(z_1,z_2,z_3,z_4)=5\pi i-\frac{1}{4}(\pi^2+z_1^2)+z_3-2z_4+\frac{1}{4}(z_1-\pi i)(z_2-z_1+z_3+z_4)+(z_1-\pi i)e^{3(z_2-z_1)+5z_3+z_4+7}-\left(e^{3(z_2-z_1)+5z_3+z_4+7}+\frac{1}{4}(z_2-z_1+z_3+z_4-2\pi i)\right)^2$ is a transcendental entire function on $\mathbb{C}^4$ with order $1$ and satisfies the 
equation 
\beas \left(\frac{\partial f(z_1,\ldots,z_4)}{\partial z_1}+\frac{\partial f(z_1,\ldots,z_4)}{\partial z_2}\right)^2+f(z_1+c_1,\ldots,z_4+c_4)=z_3-2z_4,\eeas
where $c=(\pi i,2\pi i,-\pi i, 2\pi i)$. \end{ex}
\begin{ex} Let us consider a transcendental entire function on $\mathbb{C}^5$ such that 
\beas&& f(z_1,z_2,\ldots,z_5)=\frac{4+\pi^2}{4}-\frac{1}{4}z_1^2+(z_1+\pi)\sin\left(i(z_2-z_1)+z_3+z_4-z_5\right)\\
&&-\frac{1+i}{8}(z_1+\pi)(z_2-z_1+z_3+z_4+z_5)-\left[\sin\left(i(z_2-z_1)+z_3+z_4-z_5\right)\right.\\&&\left.-\frac{1+i}{8}\left(z_2-z_1+z_3+z_4+z_5-2\pi(1-i)\right)\right]^2.\eeas 
It is easy to see that $\rho(f)=1$ and $f$ satisfies the equation 
\beas \left(\frac{\partial f(z_1,z_2,\ldots,z_5)}{\partial z_1}+\frac{\partial f(z_1,z_2,\ldots,z_5)}{\partial z_2}\right)^2+f(z_1+c_1,z_2+c_2,\ldots,z_5+c_5)=1,\eeas
where $c=(-\pi,\pi,-2\pi i,\pi, -\pi )$.\end{ex}
It is clear that we have solved the \textrm{Open problem \ref{op1}} in \textrm{Corollaries} \ref{cor1} and \ref{cor2}.
\begin{remark} The key tools in the proof of main theorems are the core part of Nevanlinna's theory, 
the difference analogue of the lemma on the logarithmic derivative in several complex variables \cite{4,14} and the Lagrange's auxiliary equations \cite[Chapter 2]{555} for quasi-linear partial differential equations.\end{remark}
\section { Some lemmas} 
The following are relevant lemmas of this paper and are used in the sequel.
\begin{lem}\label{lem5}\cite{4,14} Let $f$ be a non-constant meromorphic function with finite order on $\mathbb{C}^n$ such that $f(0)\not=0,\infty$. Then for $c\in\mathbb{C}^n$, 
\beas m\left(r, \frac{f(z)}{f(z+c)}\right)+m\left(r, \frac{f(z+c)}{f(z)}\right)=S(r, f)\eeas
holds for all $r>0$ outside of a possible exceptional set $E\subset [1,\infty)$ of finite logarithmic measure $\int_E dt/t<+\infty$.\end{lem}
\begin{lem}\label{lem6}\cite{3,25} Let $f$ be a non-constant meromorphic function with finite order on $\mathbb{C}^n$ and $I=(i_1,i_2,\ldots,i_n)$ be a multi-index with length $\Vert I\Vert=\sum_{j=1}^n i_j$. Assume
that $T(r_0, f)\geq e$ for some $r_0$. Then $m\left(r, \frac{\pa^I f}{f}\right)=S(r,f)$ holds for all $r\geq r_0$ outside a set $E\subset (0,\infty)$ of finite logarithmic measure $\int_E dt/t<+\infty$,
 where  $\pa^I f=\frac{\pa^{\Vert I\Vert} f}{\pa z_1^{i_1}\cdots\pa z_n^{i_n}}$.\end{lem}
\begin{lem}\label{lem10}\cite[\textrm{Lemma 5.34}]{12} Let $f(z)$ be a $\nu$-valued algebroid solution of the following partial differential equation 
\beas \Omega\left(z,f,\pa^{\alpha_1}f,\ldots,\pa^{\alpha_n}f\right)=\frac{A(z,f)}{B(z,f)},\eeas
where $\Omega(z,f,\pa^{\alpha_1}f,\ldots,\pa^{\alpha_n}f)=\sum_{i\in I}c_i(z)f^{i_0}\left(\pa^{\alpha_1}f\right)^{i_1}\cdots\left(\pa^{\alpha_n}f\right)^{i_n}$ with $I=\{i=\\(i_0,i_1,\ldots,i_n)\}$ is a 
finite set of distinct elements in $\mathbb{Z}_+^{n+1}$ and $c_i\in\mathscr{M}(\mathbb{C}^n)$, and $A(z,f)$ and $B(z,f)$ are co prime polynomials for $f$ given by $A(z,f)=\sum_{j=0}^p a_j(z)f^j$, 
$B(z,f)=\sum_{k=0}^qb_k(z)f^k$, where $a_j,b_k\in\mathscr{M}(\mathbb{C}^n)$ such that $a_p\not\equiv 0, b_q\not\equiv 0$. If $q\geq p$, then 
\begin{small}\beas m(r, \Omega)=O\left\{\sum_{i\in I}m\left(r,c_i\right)+\sum_{j=0}^p m\left(r,a_j\right)+\sum_{k=0}^q m\left(r,b_k\right)+m\left(r, \frac{1}{b_q}\right)+\sum_{k=1}^nm\left(r, \frac{\pa^{\alpha_k}f}{f}\right)\right\}.\eeas\end{small}
\end{lem}
\begin{lem}\label{lem20}\cite[\textrm{Lemma 3.2}, P. 385]{51} Let $f$ be a non-constant meromorphic function on $\mathbb{C}^n$. Then for any $I\in \mathbb{Z}_+^n$, $T(r, \pa^I f)=O(T(r,f))$ 
for all $r$ except possibly a set of finite Lebesgue measure, where $I=\left(i_1,i_2,\ldots,i_n\right)\in\mathbb{Z}_+^n$ denotes a multiple index with $\Vert I\Vert=i_1+i_2+\cdots+i_n$, $\mathbb{Z}_+=\mathbb{N}\cup\{0\}$, and $\pa^I f=\frac{\pa^{\Vert I\Vert} f}{\pa z_1^{i_1}\cdots\pa z_n^{i_n}}$.  
\end{lem}
\section{Proofs of the main results}  
\begin{proof}[\bf{Proof of Theorem \ref{th}}] 
Let $f$ be a finite order transcendental meromorphic function on $\mathbb{C}^n$ with $N(r,f)=S(r,f)$ satisfies (\ref{fg}) and $G(z)$ be defined in (\ref{re1}). 
Then $G(z), f(z+c)$ are finite order transcendental meromorphic functions with $N(r,G(z))=S(r,f)=N(r,f(z+c))$.
 In view of \textrm{Lemma \ref{lem5}}, we deduce that
\bs\beas\begin{array}{lll} T(r, f)=m(r, f)+S(r,f)\leq m\left(r, \frac{f(z)}{ f(z+c)}\right)+m(r, f(z+c))+S(r,f)\leq T(r, f(z+c))+S(r,f),\nonumber\\
T(r, f(z+c))=m(r, f(z+c))+S(r,f)\leq m\left(r, \frac{f(z+c)}{ f(z)}\right)+m(r, f)+S(r,f)\leq T(r, f)+S(r,f).\end{array}\eeas\es
Thus, $T(r,f(z+c))=T(r,f)+S(r,f)$.
Note that, $m\left(r,\sum\limits_{m=1}^{n}\sum\limits_{\Vert I\Vert=m}a_{I}(z)\frac{\pa^I f(z)}{f(z)}\right)=S(r,f)$, since $f$ is a finite order transcendental meromorphic function and $a_I(z)$ are small functions of $f$, where $I=(i_1,\ldots,i_n)\in\mathbb{Z}_+^{n}$ with $\Vert I\Vert=\sum_{j=1}^n i_j$. Then by \textrm{Lemma \ref{lem5}} and \ref{lem20}, we have
\bea\label{ree4} T(r,G(z))=m(r,G(z))+S(r,f)&=&m\left(r,f(z)\sum\limits_{m=1}^{n}\sum\limits_{|\Vert I\Vert=m}a_{I}(z)\frac{\pa^I f(z)}{f(z)}\right)+S(r,f)\nonumber\\
&\leq& T(r,f(z))+S(r,f).\eea
Now we discuss the following two cases.\\
{\bf Case 1.}  When $m_2>m_1$.
In view of Valiron-Mokhon'ko lemma \cite[p. 29]{12} and \textrm{Lemma \ref{lem6}} and (\ref{ree4}), we have 
\beas m_2T(r, f(z))&=& m_2T(r, f(z+c))+S(r,f)\\
&=&T(r, f^{m_2}(z+c))+S(r,f) \\
&\leq& T\left(r, \alpha(z)f^{m_2}(z+c)\right)+S(r,f)\\
&=&T\left(r, G^{m_1}(z)-\beta(z) \right)+S(r,f)\\
&\leq&m_1T\left(r,G)\right)+S(r,f)\\
&\leq& m_1T(r, f(z))+S(r,f),\\\text{\it i.e.,}\quad
(m_2-m_1)T(r, f(z))&\leq& S(r,f),\eeas which arise a contradiction, since $f$is a finite order transcendental meromorphic function and $m_2>m_1$.\\
{\bf Case 2.} When $m_1>m_2\geq 2$ {\it i.e.,} $\frac{1}{m_1}+\frac{1}{m_2}<1$, which implies $m_2>\frac{m_1}{m_1-1}$.
By Nevanlinna second fundamental theorem for small functions \cite[p. 50]{12}, \textrm{Lemma \ref{lem5}} and (\ref{fg}), we have
\bea\label{fg2} m_1T(r, G(z))&=&T\left(r, G^{m_1}(z)\right)+S(r,G)\nonumber\\
&\leq&\ol N\left(r, G^{m_1}(z)\right)+\ol N\left(r, 0;G^{m_1}(z)\right)\nonumber\\
&&+\ol N\left(r,0;G^{m_1}(z)-\beta(z)\right)+S(r,f)\nonumber\\
&=&\ol N\left(r, 0;G(z)\right)+\ol N\left(r,0; \alpha(z)f^{m_2}(z+c)\right)+S(r,f)\nonumber\\
&\leq& T(r, G(z))+\ol N\left(r, 0; f(z+c)\right)+S(r,f),\nonumber\\\text{\it i.e.,}\;\;
(m_1-1)T(r, G(z))&\leq& T(r, f(z+c))+S(r,f).\eea
Again, in view of Valiron-Mokhon'ko lemma \cite[p. 29]{12}, \textrm{Lemma \ref{lem5}}, (\ref{fg}) and (\ref{fg2}), we have
\beas m_2T(r, f(z+c))&=& T(r, \alpha(z)f^{m_2}(z+c))+S(r,f)\\
&=&T(r, G^{m_1}(z)-\beta(z))+S(r,f)\\
&\leq&m_1T(r, G(z))+S(r,f)\\
&\leq&\frac{m_1}{m_1-1}T(r, f(z+c))+S(r,f),\\\text{\it i.e.,}\quad
\left(m_2-\frac{m_1}{m_1-1}\right)T(r, f(z+c))&\leq&S(r,f),\eeas 
which arise a contradiction, since $f$ is a finite order transcendental meromorphic function and $m_2>\frac{m_1}{m_1-1}$. This completes the proof. 
\end{proof}
\begin{proof}[\bf{Proof of Theorem \ref{T1}}] Let $f$ be a finite order transcendental entire function on $\mathbb{C}^n$ satisfies (\ref{fte}). 
Differentiating partially with respect to $z_1$ on both sides of (\ref{fte}), we have 
\bea\label{2r} m_1\left(\frac{\pa f(z)}{\pa z_1}\right)^{m_1-1}\frac{\pa^2 f(z)}{\pa z_1^2}+\frac{\pa f(z+c)}{\pa z_1}-\frac{\pa f(z)}{\pa z_1}=0.\eea
Let $F(z)=\frac{\pa f(z)}{\pa z_1}$. Then (\ref{2r}) reduces to 
\bea\label{fte3}m_1F^{m_1-1}(z)\frac{\pa F(z)}{\pa z_1}=-F(z+c)+F(z),\nonumber\\
\text{\it{i.e.,}}\quad F^{m_1-2}(z)\frac{\pa F(z)}{\pa z_1}=-\frac{1}{m_1}\frac{F(z+c)-F(z)}{F(z)}.\eea
Clearly by \textrm{Lemmas \ref{lem5}} and \ref{lem20}, we get $m\left(r, -\frac{1}{m_1}\frac{F(z+c)-F(z)}{F(z)}\right)=S(r, F)$, which implies
 $m\left(r, F^{m_1-2}(z)\frac{\pa F(z)}{\pa z_1}\right)=S(r, F)=S(r, f)$.
Since $f$ is a finite order transcendental entire function on $\mathbb{C}^n$, we see that $N\left(r, F^{m_1-2}(z)\frac{\pa F(z)}{\pa z_1}\right)=S(r, f)$. Consequently, we get $T\left(r, F^{m_1-2}(z)\frac{\pa F(z)}{\pa z_1}\right)=S(r, f)$.
Since $f$ is a finite order transcendental entire function and from (\ref{fte3}), we may assume that 
\bea\label{fte1}F^{m_1-2}(z)\frac{\pa F(z)}{\pa z_1}\equiv P(z),\eea 
where $P(z)$ is a non-zero polynomial on $\mathbb{C}^n$. 
Then the following two cases arise.\\
{\bf Case 1.} When $m_1=2$. Solving (\ref{fte1}) by using the Lagrange method \cite[Chapter 2]{555}, we get $F(z)=Q(z)+g_1(z_2,z_3,\ldots,z_n)$, 
where $g_1(z_2,z_3,\ldots,z_n)$ is a finite order entire function in $z_2,z_3,\ldots,z_n$ and $Q(z)=\int P(z)dz_1$, where $z_2,z_3,\ldots,z_n$ are constants. Note that $\deg (Q(z))\geq 1$.
Now the following two cases arise.\\
{\bf Sub-case 1.1.} Let $g_1(z_2,z_3,\ldots,z_n)$ be a finite order transcendental entire function in $z_2,z_3,\ldots,z_n$. From (\ref{fte3}), we have
\bea\label{fte4}&&\frac{F(z+c)}{F(z)}\equiv 1-2P(z),\;\;\text{\it i.e.,}\;\;\frac{Q(z+c)+g_1(z_2+c_2,\ldots,z_n+c_n)}{Q(z)+g_1(z_2,z_3,\ldots,z_n)}\equiv1-2P(z)\nonumber\\[2mm]\text{\it i.e.,}
&&Q(z+c)+g_1(z_2+c_2,\ldots,z_n+c_n)\equiv(1-2P(z))Q(z)\nonumber\\[2mm]&&+(1-2P(z))g_1(z_2,z_3,\ldots,z_n).\eea
Comparing the polynomials on the both sides, we get $P(z)\equiv -1/2$ and $Q(z)=-z_1/2$. 
Hence, we have  
\bea\label{R2}&& F(z)=-\frac{1}{2}z_1+g_1(z_2,z_3,\ldots,z_n)\\\text{and}
\label{R1}&&g_1(z_2+c_2,z_3+c_3,\ldots,z_n+c_n)\equiv g_1(z_2,z_3,\ldots,z_n)+\frac{c_1}{2}.\eea
From (\ref{R1}), we deduce that $g_1(z_2,z_3,\ldots,z_n)\equiv g_2(z_2,z_3,\ldots,z_n)+c_1\omega/(2\tau)$, where $g_2$ is a finite order transcendental entire periodic function in $z_2,z_3,\ldots,z_n$ with period $(c_2,c_3,\ldots,c_n)\in\mathbb{C}^{n-1}\setminus\{0\}$, $\tau=c_2+c_3+\cdots+c_n\not=0$ and $\omega=z_2+z_3+\cdots+z_n$.
On integration from (\ref{R2}), we have
\bea\label{fte5} f(z)=-\frac{1}{4}z_1^2+z_1g_1(z_2,z_3,\ldots,z_n)+g_3(z_2,z_3,\ldots,z_n),\eea 
where $g_3(z_2,z_3,\ldots,z_n)$ is a finite order entire function in $z_2,z_3,\ldots,z_n$.\\
Using $g_1(z_2,z_3,\ldots,z_n)\equiv g_2(z_2,z_3,\ldots,z_n)+\frac{c_1}{2\tau}\omega$, we deduce from (\ref{fte}) and (\ref{fte5}) that
\beas &&\left(-\frac{1}{2}z_1+g_1(z_2,z_3,\ldots,z_n)\right)^2-\frac{1}{4}(z_1+c_1)^2+(z_1+c_1)g_1(z_2+c_2,z_3+c_3,\ldots,\\&&z_n+c_n)+g_3(z_2+c_2,z_3+c_3,\ldots,z_n+c_n)\equiv \varphi(z_2,z_3,\ldots,z_n),\\
\text{\it i.e.,} &&g_3(z_2,z_3,\ldots,z_n)\equiv \varphi(z_2-c_2,z_3-c_3,\ldots,z_n-c_n)+\frac{c_1^2}{4}\\&&-\left(g_2(z_2,z_3,\ldots,z_n)+\frac{c_1}{2\tau}(\omega-\tau)\right)^2-c_1\left(g_2(z_2,z_3,\ldots,z_n)+\frac{c_1\omega}{2\tau}\right).\eeas 
Thus,
\beas &&f(z)=(z_1-c_1)\left(g_2(z_2,z_3,\ldots,z_n)+\frac{c_1}{2\tau}\omega\right)-\left(g_2(z_2,z_3,\ldots,z_n)+\frac{c_1}{2\tau}(\omega-\tau)\right)^2\\&&+\frac{1}{4}\left(c_1^2-z_1^2\right)+\varphi(z_2-c_2,z_3-c_3,\ldots,z_n-c_n),\eeas 
where $g_2(z_2,z_3,\ldots,z_n)$ is a finite order transcendental entire periodic function in $z_2,z_3,\\\ldots,z_n$ with period $(c_2,c_3,\ldots,c_n)\in\mathbb{C}^{n-1}$, $\varphi(z_2,z_3,\ldots,z_n)$ is a finite order entire function, $\omega=z_2+z_3+\cdots+z_n$ and $\tau=c_2+c_3+\cdots+c_n\not=0$.\\
{\bf Sub-case 1.2.} Let $g_1(z_2,z_3,\ldots,z_n)$ be a polynomial in $z_2,z_3,\ldots,z_n$. 
Now proceeding similarly as \textrm{Sub-case 1.1}, we again get $f(z)=-\frac{1}{4}z_1^2+z_1g_1(z_2,z_3,\ldots,z_n)+g_3(z_2,z_3,\ldots,z_n)$, where $g_3(z_2,z_3,\ldots,z_n)$ is a finite order transcendental entire function in $z_2,z_3,\ldots,z_n$ and $g_1$ satisfies (\ref{R1}).
From (\ref{fte}), we deduce that
\bea\label{fte7} &&\left(-\frac{1}{2}z_1+g_1(z_2,z_3,\ldots,z_n)\right)^2-\frac{1}{4}(z_1+c_1)^2+(z_1+c_1)g_1(z_2+c_2,z_3+c_3,\ldots,\nonumber\\&&z_n+c_n)+g_3(z_2+c_2,z_3+c_3,\ldots,z_n+c_n)\equiv \varphi(z_2,z_3,\ldots,z_n),\nonumber\\
\text{\it i.e.,} &&g_3(z_2+c_2,\ldots,z_n+c_n)\equiv\varphi(z_2,z_3,\ldots,z_n) -\left(-\frac{1}{2}z_1+g_1(z_2,z_3,\ldots,z_n)\right)^2\nonumber\\
&&+\frac{1}{4}(z_1+c_1)^2-(z_1+c_1)g_1(z_2+c_2,z_3+c_3,\ldots,z_n+c_n),\nonumber\\
\text{\it i.e.,} &&g_3(z_2,\ldots,z_n)\equiv\varphi(z_2-c_2,z_3-c_3,\ldots,z_n-c_n)+\frac{1}{4}z_1^2-z_1g_1(z_2,z_3,\ldots,z_n)\nonumber\\
&& -\left(-\frac{1}{2}(z_1-c_1)+g_1(z_2-c_2,z_3-c_3,\ldots,z_n-c_3)\right)^2.\eea 
Since $g_1(z_2,z_3,\ldots,z_n)$ is a polynomial in $z_2,z_3,\ldots,z_n$ while $g_3(z_2,z_3,\ldots,z_n)$ is a finite order transcendental entire function in $z_2,z_3,\ldots,z_n$. From (\ref{fte7}), we must have $\varphi(z_2,z_3,\ldots,z_n)$ is a finite order transcendental entire function, otherwise contradiction arise. Thus
\bs\beas f(z)=\varphi(z_2-c_2,z_3-c_3,\ldots,z_n-c_n)-\left(-\frac{1}{2}(z_1-c_1)+g_1(z_2-c_2,z_3-c_3,\ldots,z_n-c_3)\right)^2,\eeas \es
where $g_1(z_2,z_3,\ldots,z_n)$ is a polynomial such that $g_1(z_2+c_2,z_3+c_3,\ldots,z_n+c_n)\equiv g_1(z_2,z_3,\ldots,z_n)+\frac{c_1}{2}$ holds and $\varphi(z_2,z_3,\ldots,z_n)$ is a finite order transcendental entire function in $z_2,z_3,\ldots,z_n$.\\  
{\bf Case 2.} When $m_1\geq 3$. If $m\left(r,P(z)\right)\not=S(r,F)$, i.e., $T\left(r,P(z)\right)\not=S(r,F)$, then we conclude from (\ref{fte1}) that $F(z)$ and $\frac{\pa F(z)}{\pa z_1}$ are both non-zero polynomials on $\mathbb{C}^n$. Otherwise, if $F(z)$ or $\frac{\pa F(z)}{\pa z_1}$ is transcendental, then L.H.S. of (\ref{fte1}) is transcendental while its R.H.S. is polynomial and it is not possible. 
If $m\left(r,P(z)\right)=S(r,F)$, by \textrm{Lemma \ref{lem10}} and (\ref{fte1}), we have $m\left(r,\frac{\pa F(z)}{\pa z_1}\right)=S(r, F)=S(r, f)$. Since $f$ is a finite order transcendental entire function on $\mathbb{C}^n$, we see that $T\left(r,\frac{\pa F(z)}{\pa z_1}\right)=S(r, f)$ and so in view of (\ref{fte1}), we may assume that 
\bea\label{fte6}\frac{\pa F(z)}{\pa z_1}\equiv Q(z)\quad\text{which implies}\quad F(z)=R(z)+g_4(z_2,z_3,\ldots,z_n),\eea 
where $Q(z)$ is a non-zero polynomial on $\mathbb{C}^n$, $R(z)=\int Q(z)dz_1$, in which $z_2,z_3,\ldots,z_n$ are constants with $\deg\left(R(z)\right)\geq 1$ and $g_4(z_2,z_3,\ldots,z_n)$ is a finite order entire function in $z_2,z_3,\ldots,z_n$. Now the following cases arise.\\
{\bf Sub-case 2.1.} Let $g_4(z_2,z_3,\ldots,z_n)$ be a finite order transcendental entire function in $z_2,z_3,\ldots,z_n$.
Since $m_1\geq 3$, we see that L.H.S. of (\ref{fte1}) is a finite order transcendental entire function while its R.H.S. is a polynomial and this arise a contradiction.\\ 
{\bf Sub-case 2.2.} Let $g_4(z_2,z_3,\ldots,z_n)$ be a polynomial in $z_2,z_3,\ldots,z_n$. Then $F(z)$ is a polynomial with $\deg(F(z))\geq 1$.
From (\ref{fte3}) and (\ref{fte1}), we have
\beas
R(z+c)+g_4(z_2+c_2,z_3+c_3,\ldots,z_n+c_n)&\equiv& -m_1P(z)\left(R(z)+g_4(z_2,z_3,\ldots,z_n)\right).\eeas
Comparing the polynomials on the both sides, we get $P(z)\equiv -1/m_1$. Then from (\ref{fte1}), we see that
\beas\label{fte8} \left(F(z)\right)^{m_1-2}Q(z)\equiv -\frac{1}{m_1},\eeas
which arise a contradiction by comparing the degrees on the both sides. This completes the proof.
\end{proof}
\begin{proof}[\bf{Proof of Theorem \ref{T2}}] 
Let $f$ be a finite order transcendental entire function on $\mathbb{C}^n$ satisfies (\ref{ftee}) and $F(z)=\left(\frac{\pa }{\pa z_1}+\frac{\pa }{\pa z_2}\right)f(z)$. Differentiating partially with respect to $z_1$ and $z_2$ respectively on both sides of (\ref{ftee}), we have
\bea\label{2er1} m_1F^{m_1-1}(z)\frac{\pa F(z)}{\pa z_1}=-\frac{\pa f(z+c)}{\pa z_1}\;\;\text{and}\;\;m_1F^{m_1-1}(z)\frac{\pa F(z)}{\pa z_2}=-\frac{\pa f(z+c)}{\pa z_2}\eea
From (\ref{2er1}), we have 
\bea\label{2er2} m_1F^{m_1-1}(z)\left(\frac{\pa }{\pa z_1}+\frac{\pa }{\pa z_2}\right)F(z)=-F(z+c)
.\eea
By \textrm{Lemmas \ref{lem5}} and \ref{lem20}, and from  (\ref{2er2}), we get 
\beas m\left(r, F^{m_1-2}(z)\left(\frac{\pa F(z)}{\pa z_1}+\frac{\pa F(z)}{\pa z_2}\right)\right)=m\left(r, -\frac{1}{m_1}\frac{F(z+c)}{F(z)}\right)=S(r, F)=S(r, f).\eeas
Note that $T\left(r, F^{m_1-2}(z)\left(\frac{\pa F(z)}{\pa z_1}+\frac{\pa F(z)}{\pa z_2}\right)\right)=S(r, f)$, since $f$ is a finite order transcendental entire function and from (\ref{2er2}), we may assume that 
\bea\label{2fte1}F^{m_1-2}(z)\left(\frac{\pa F(z)}{\pa z_1}+\frac{\pa F(z)}{\pa z_2}\right)\equiv P(z),\eea 
where $P(z)$ is a non-zero polynomial on $\mathbb{C}^n$.
Now the following two cases arise.\\[2mm]
{\bf Case 1.} When $m_1=2$. 
Then from (\ref{2fte1}), we have
\bea\label{wr5}\frac{\pa F(z)}{\pa z_1}+\frac{\pa F(z)}{\pa z_2}\equiv P(z), \eea
where $P(z)$ is a non-zero polynomial on $\mathbb{C}^n$. The Lagrange's auxiliary equations \cite[Chapter 2]{555} corresponding to (\ref{wr5}) are as follows
\beas\frac{dz_1}{1}=\frac{dz_2}{1}=\frac{dz_3}{0}=\cdots=\frac{dz_n}{0}=\frac{dF}{P(z)}.\eeas
Note that $a_2=z_2-z_1$, $a_i=z_i$ ($3\leq i\leq n$) and $dF=P(z)dz_1=P(z_1,z_1+a_2,a_3,\cdots,a_n)dz_1$ which implies $F(z)=Q(z)+a_1$, where $Q(z)$ is obtained replacing $a_2$ by $z_2-z_1$, $a_3$ by $z_3$, ... , $a_n$ by $z_n$ in the integration of $P(z_1,z_1+a_2,a_3,\cdots,a_n)$ w.r.t. $z_1$ and $a_i\in\mathbb{C}$ ($1\leq i\leq n$). 
Hence the solution is $\phi(a_1,\cdots,a_n)=0$. For simplicity, we suppose 
\beas F(z)=Q(z)+g_2(z_2-z_1,z_3,\cdots,z_n), \eeas
where $g_2(z_2-z_1,z_3,\cdots,z_n)$ is a finite order entire function in $z_2-z_1,z_3,\cdots,z_n$ and $Q(z)$ is a non-zero polynomial on $\mathbb{C}^n$ with $\deg (Q(z))\geq 1$. Hence, we obtain 
\bea\label{twr6}\frac{\pa f(z)}{\pa z_1}+\frac{\pa f(z)}{\pa z_2}=Q(z)+g_2(z_2-z_1,z_3,\cdots,z_n).\eea
The Lagrange's auxiliary equations \cite[Chapter 2]{555} corresponding to (\ref{twr6}) are as follows 
\bea\label{twr8} \frac{dz_1}{1}=\frac{dz_2}{1}=\frac{dz_3}{0}=\cdots=\frac{dz_n}{0}=\frac{df}{Q(z)+g_2(z_2-z_1,z_3,\cdots,z_n)}.\eea
Now the following two cases arise.\\
{\bf Sub-case 1.1.} Let $g_2(z_2-z_1,z_3,\cdots,z_n)$ be a polynomial in $z_2-z_1,z_3,\ldots,z_n$.
Then from (\ref{twr8}), we have $z_2-z_1=d_2$, $z_i=d_i$ ($3\leq i\leq n$) and $df=Q(z_1,z_1+d_2,d_3,\cdots,d_n)dz_1+g_2(d_2,d_3,\cdots,d_n)dz_1$ which implies \beas f(z)=R(z)+z_1g_2(z_2-z_1,z_3,\cdots,z_n)+d_1,\eeas 
where $R(z)$ is obtained replacing $d_2$ by $z_2-z_1$, $d_3$ by $z_3$, $\cdots$ , $d_n$ by $z_n$ in the integration of $Q(z_1,z_1+d_2,d_3,\cdots,d_n)$ w.r.t. $z_1$ and $d_i\in\mathbb{C}$ ($1\leq i\leq n$). Hence the solution is $\psi(d_1,\cdots,d_n)=0$. For simplicity, we suppose 
\beas f(z)=R(z)+z_1g_2(z_2-z_1,z_3,\cdots,z_n)+g_3(z_2-z_1,z_3,\cdots,z_n),\eeas 
where $g_3(z_2-z_1,z_3,\cdots,z_n)$ is a finite order entire function in $z_2-z_1,z_3,\cdots,z_n$. Since $f$ is transcendental entire function, so we must have $g_3(z_2-z_1,z_3,\cdots,z_n)$ is a finite order transcendental entire function. From (\ref{ftee}), we have 
\bs\bea &&\left(Q(z)+g_2(z_2-z_1,z_3,\cdots,z_n)\right)^2+R(z+c)+(z_1+c_1)g_2(z_2-z_1+c_2-c_1,\nonumber\\
&&z_3+c_3,\cdots,z_n+c_n)+g_3(z_2-z_1+c_2-c_1,z_3+c_3,\cdots,z_n+c_n)\equiv \varphi(z_3,z_4,\ldots,z_n),\nonumber\\\text{\it i.e.,}
&&g_3(z_2-z_1,z_3,\cdots,z_n)\equiv \varphi(z_3-c_3,z_4-c_4,\ldots,z_n-c_n)-z_1g_2(z_2-z_1,z_3,\cdots,z_n)\nonumber\\
\label{R3}&&-R(z)-\left(Q(z-c)+g_2(z_2-z_1-c_2+c_1,z_3-c_3,\cdots,z_n-c_n)\right)^2.\eea\es
Since $g_3$ is a finite order transcendental entire function, so we must have $\varphi(z_3,z_4,\ldots,z_n)$ is a finite order transcendental entire function, otherwise we get a contradiction from (\ref{R3}). From (\ref{2er2}), we deduce that
\bs\be\label{R4} 2\left(Q(z)+g_2(z_2-z_1,z_3,\cdots,z_n)\right)P(z)\equiv -Q(z+c)-g_2(z_2-z_1+c_2-c_1,z_3+c_3,\cdots,z_n+c_n).\ee\es
From (\ref{R4}), we have $P(z)\equiv-1/2$ and hence $Q(z)=-z_1/2$, $R(z)=-z_1^2/4$ and 
\bea\label{R5} g_2(z_2-z_1+c_2-c_1,z_3+c_3,\cdots,z_n+c_n)\equiv g_2(z_2-z_1,z_3,\cdots,z_n)+\frac{c_1}{2}.\eea
Note that, from (\ref{wr5}), we have $\frac{\pa g_2}{\pa z_1}+\frac{\pa g_2}{\pa z_2}\equiv 0$.
Thus
\beas &&f(z)=-\left(-\frac{1}{2}(z_1-c_1)+g_2(z_2-z_1-c_2+c_1,z_3-c_3,\cdots,z_n-c_n)\right)^2\\&&+\varphi(z_3-c_3,z_4-c_4,\ldots,z_n-c_n),\eeas
where $g_2$ is a polynomial satisfies (\ref{R5}) with $\frac{\pa g_2}{\pa z_1}+\frac{\pa g_2}{\pa z_2}\equiv 0$ and $\varphi(z_3,z_4,\ldots,z_n)$ is a finite order transcendental entire function.\\
{\bf Sub-case 1.2.} Let $g_2(z_2-z_1,z_3,\cdots,z_n)$ be a finite order transcendental entire function in $z_2-z_1,z_3,\ldots,z_n$. 
Similarly as \textrm{Sub-case 1.1.}, we deduce that $f(z)=R(z)+z_1g_2(z_2-z_1,z_3,\cdots,z_n)+g_3(z_2-z_1,z_3,\cdots,z_n)$, where $g_3(z_2-z_1,z_3,\cdots,z_n)$ is a finite order entire function in $z_2-z_1,z_3,\cdots,z_n$.
Similarly, from (\ref{2er2}), we obtain the equation (\ref{R4}) and hence 
we have $P(z)\equiv-1/2$ and hence $Q(z)=-z_1/2$ and $R(z)=-z_1^2/4$. Therefore $g_2(z_2-z_1+c_2-c_1,z_3+c_3,\cdots,z_n+c_n)\equiv g_2(z_2-z_1,z_3,\cdots,z_n)+\frac{c_1}{2}$. 
Now we deduce that $g_2(z_2-z_1,z_3,\cdots,z_n)\equiv g_4(z_2-z_1,z_3,\cdots,z_n)+\frac{c_1}{2\tau}(z_2-z_1+z_3+\cdots+z_n)$, where $g_4$ is a finite order transcendental entire periodic function in $z_2-z_1,z_3,\ldots,z_n$ with period $(c_2-c_1,c_3,\ldots,c_n)\in\mathbb{C}^{n-1}$ and $\tau=c_2-c_1+c_3+c_4+\cdots+c_n\not=0$. 
Note from (\ref{wr5}) that $\frac{\pa g_4}{\pa z_1}+\frac{\pa g_4}{\pa z_2}\equiv 0$.
Therefore, we have
\bea\label{twr81} &&f(z)=-\frac{1}{4}z_1^2+z_1\left(g_4(z_2-z_1,z_3,\cdots,z_n)+\frac{c_1}{2\tau}(z_2-z_1+z_3+\cdots+z_n)\right)\nonumber\\&&+g_3(z_2-z_1,z_3,\cdots,z_n),\eea
where $g_3$ is a finite order entire function in $z_2-z_1,z_3,\ldots,z_n$. Putting (\ref{twr81}) into (\ref{ftee}), we get
\beas
&&g_3(z_2-z_1,z_3,\cdots,z_n)\equiv \varphi(z_3-c_3,z_4-c_4,\ldots,z_n-c_n)-\frac{1}{4}c_1^2\\
&&-c_1\left(g_4(z_2-z_1,z_3,\cdots,z_n)+\frac{c_1}{2\tau}(z_2-z_1+z_3+\cdots+z_n-\tau)\right)\\&&-\left(g_4(z_2-z_1,z_3,\cdots,z_n)+\frac{c_1}{2\tau}(z_2-z_1+z_3+\cdots+z_n-\tau)\right)^2.\eeas
Therefore, we get from (\ref{twr81}) that
\beas &&f(z)=\frac{1}{4}\left(c_1^2-z_1^2\right)+(z_1-c_1)\left(g_4(z_2-z_1,z_3,\cdots,z_n)+\frac{c_1}{2\tau}(z_2-z_1+z_3+\cdots+z_n)\right)\\&&-\left(g_4(z_2-z_1,z_3,\cdots,z_n)+\frac{c_1}{2\tau}(z_2-z_1+z_3+\cdots+z_n-\tau)\right)^2\\&&+\varphi(z_3-c_3,z_4-c_4,\ldots,z_n-c_n),\eeas
where $g_4(z_2-z_1,z_3,\cdots,z_n)$ is a finite order transcendental entire function with period \\$(c_2-c_1,c_3,\ldots,c_n)$ with $\frac{\pa g_4}{\pa z_1}+\frac{\pa g_4}{\pa z_2}\equiv0$ and $\varphi(z_3,z_4,\ldots,z_n)$ is a finite order entire function.\\
{\bf Case 2.} When $m_1\geq 3$. If $m\left(r,P(z)\right)\not=S(r,F)$, i.e., $T\left(r,P(z)\right)\not=S(r,F)$, then we conclude from (\ref{2fte1}) that $F(z)$ and $\left(\frac{\pa }{\pa z_1}+\frac{\pa }{\pa z_2}\right)F(z)$ are both non-zero polynomials on $\mathbb{C}^n$. Otherwise, if $F(z)$ or $\left(\frac{\pa }{\pa z_1}+\frac{\pa }{\pa z_2}\right)F(z)$ 
is a finite order transcendental entire function in $z_1,z_2$ and other variables, we see that L.H.S. of (\ref{2fte1}) is a finite order transcendental entire function while its R.H.S. is a polynomial and this is not possible. 
If $m\left(r,P(z)\right)=S(r,F)$, by \textrm{Lemma \ref{lem10}} and (\ref{2fte1}), we have $m\left(r,\frac{\pa F(z)}{\pa z_1}+\frac{\pa F(z)}{\pa z_2}\right)=S(r, f)$.
Since $f$ is a finite order transcendental entire function, we have $T\left(r,\frac{\pa F(z)}{\pa z_1}+\frac{\pa F(z)}{\pa z_2}\right)=S(r, f)$ and from (\ref{2fte1}), we may assume that 
\beas\frac{\pa f(z)}{\pa z_1}+\frac{\pa f(z)}{\pa z_2}\equiv P_1(z)\;\text{and}\;\frac{\pa F(z)}{\pa z_1}+\frac{\pa F(z)}{\pa z_2}\equiv P_2(z),\eeas 
where $P_1(z)$ and $P_2(z)$ are non-zero polynomials on $\mathbb{C}^n$. The remaining part of the proof follows from \textrm{Sub-case 1.1.} of this theorem. 
This completes the proof.\end{proof}
\section{Declarations}
\noindent{\bf Acknowledgment:} 
Third Author is supported by University Grants Commission (IN) fellowship (No. F. 44 - 1/2018 (SA - III)). We would also want to thank the anonymous reviewers and the editing team for their suggestions.\\
{\bf Author's contributions:} All authors have equal contribution to complete the manuscript. All of them read and approved the final manuscript.\\
{\bf Conflict of Interest:} Authors declare that they have no conflict of interest.\\
{\bf Availability of data and materials:} Not applicable.\\

\end{document}